\documentclass[10pt,journal]{IEEEtran}

\usepackage{xspace,amssymb,epsfig,syntonly}
\usepackage{epsfig,amsmath,color}
\usepackage{xspace,syntonly,empheq}
\usepackage{wrapfig}
\usepackage{url}
\usepackage{bbm}
\usepackage{mathtools}
\usepackage{cite}
\usepackage{graphicx}
\usepackage{caption}
\usepackage{subcaption}
\usepackage{algorithm,algorithmic}

\DeclarePairedDelimiterX{\norm}[1]{\lVert}{\rVert}{#1}

\usepackage{epstopdf}

\begin{document}

\newtheorem{definition}{Definition}
\newtheorem{remark}{Remark}
\newtheorem{proposition}{Proposition}
\newtheorem{lemma}{Lemma}
\def\HS{\hspace{\fontdimen2\font}}
\font\myfont=cmr12 at 16pt

\IEEEoverridecommandlockouts

\markboth{SUBMITTED TO IEEE POWER ENGINEERING LETTERS}%
{ \MakeLowercase{\textit{et al.}}: }


\title{Directly Constraining Marginal Prices}

\author{Kyri Baker, \emph{Member, IEEE}
\thanks{K. Baker is with the National Renewable Energy Laboratory, Golden, CO. E-mail: {kyri.baker@nrel.gov}%
}
}
\maketitle

\begin{abstract}
The marginal price of electricity traditionally depends on the dual variables associated with relevant optimization goals. Particularly, in the optimal power flow realm, prices represent the cost of supplying an additional unit of power at each bus; for the economic dispatch case, dual variables represent the cost of supplying an additional unit of power to the whole system. Dual variables are useful for many additional tasks, including the analysis of system congestion and the determination of the cost of load adjustments. In this letter, departing from conventional optimal power flow analysis, it is shown how constraints on relevant dual variables affect the prices of electricity, allowing for market settings and demand response programs that accept bids and caps on individual prices. 

\end{abstract}


\section{Introduction}
\label{sec:Introduction}
As we move towards a smarter and more flexible power grid, an increased number of entities may be willing to participate in demand response (DR), from industrial customers to individual households. This heightened participation could have the potential to drastically reduce or increase the load, and perhaps even cause instability and oscillations in the system as multiple entities simultaneously respond to DR signals. This phenomenon has already been observed to some degree in current demand response programs and is dubbed the ``rebound effect" \cite{rebound}. 

In addition, many limitations exist in current DR programs \cite{duncan2011}. Many of the DR techniques involve one-way communication of a pricing signal that the customer can decide to respond to, or offer a flat rate for reducing demand during peak hours. For example, Pacific Gas \& Electric's Demand Bidding Program gives a \$0.50/kWh for providing load reduction service during specific times. In PJM Interconnection's demand bidding program, consumers make bids for a portion of their load that they are willing to reduce. The bid is accepted when the market clears if it is less than the market price. 

Current approaches may neglect to account for the fact that as an increased number of entities are participating in DR programs, the benefit from DR may decrease, the amount of volatility in the system may increase, and the ``rebound effect" may occur. For example, if a utility sends out a signal advertising a \$0.50/kWh payment for consumers to reduce their demand in the next hour and more loads than necessary decide to participate, the utility's gain from this decreases. Conversely, demand side bidding (DSB) requires each participating entity to specify a price and kWh bid, which may not be straightforward to determine, and is not always included directly into the optimization problem. For example, a homeowner may wish to spend no more than \$100 a month on electricity, but is unable to determine the corresponding level of energy to and therefore is unable to submit an informed bid in this market structure.

In this letter we propose an alternate approach to DSB that includes bids directly in the optimization problem. Issues such as price volatility and market power can be mitigated by DSB\cite{DSB}, and we propose to implement this DR structure by not only accepting (price, kWh) bids, but by accepting bids based solely on price. Mathematically, a method to include explicit constraints on electricity price has currently not been formulated due to the fact that price is usually an output of the optimization and not known a priori. However, in this letter, it can be shown that by using optimization duality theory, these prices can be directly constrained, and this new constraint translates to a new variable in terms of power consumption in the original optimization problem. When these variables can be constrained directly, demand side bidding, individual consumer budgeting, and prevention of price fluctuations can be explicitly considered in the optimization problem.

\section{Derivation of the Price-Constrained Optimization Problem}
\label{sec:QP}
In order to demonstrate how constraining the Lagrange multipliers affects the primal problem, first define a general quadratic programming problem (a general optimal power flow problem with linearized AC power flow equations or no network consideration):

\vspace{-.2cm}
\begin{small}
\begin{equation} 
\begin{aligned}
& \underset{x}{\text{minimize}}
& & \frac{1}{2}x^TQx + c^Tx \\
& \text{subject to}
& & Ax \leq b
\end{aligned}
\end{equation}
\end{small}

\noindent where $x$ is a $n$-dimensional vector of primal problem variables such as generator values, voltage angle, and flexible load values. Matrix $Q$ is an $n \times n$ real, positive definite symmetric matrix, $A$ is an $m \times n$ real matrix, and $b$ is an $m$-dimensional real vector. The linearized power flow equations can be captured in the constraints. We form the Lagrangian dual function:

\begin{small}
\begin{equation} \label{dualfcn}
\begin{aligned}
& g(\lambda) = \underset{x}{\text{inf}}
& & \frac{1}{2}x^TQx + c^Tx + \lambda^T(Ax-b)\\
\end{aligned}
\end{equation}
\end{small}
\vspace{-.2cm}

\noindent Where $\lambda$ is a $m$-dimensional vector of dual variables. We observe that the infimum is obtained for $x = -Q^{-1}(c+A^T\lambda)$. Substituting this into (\ref{dualfcn}), we obtain the dual function in terms of only the dual variables:

\begin{small}
\begin{equation} \label{dualfcn2}
\begin{aligned}
& g(\lambda) =  -\frac{1}{2}\lambda^TAQ^{-1}A^T\lambda - \lambda^T(b+AQ^{-1}c) - \frac{1}{2}c^TQ^{-1}c\\
\end{aligned}
\end{equation}
\end{small}
\vspace{-.2cm}

\noindent Defining $P = AQ^{-1}A^T$ and $t = b + AQ^{-1}c$, we can write the dual problem as the following:

\vspace{-.2cm}
\begin{small}
\begin{equation} 
\begin{aligned}
& \underset{\lambda}{\text{minimize}}
& & \frac{1}{2}\lambda^TP\lambda + t^T\lambda \\
& \text{subject to}
& & \lambda \geq 0
\end{aligned}
\end{equation}
\end{small}
\vspace{-.2cm}

\noindent  Assuming Slater's condition holds, and because the primal problem is convex, strong duality will hold between the primal and dual problems \cite{BoVa04}, and the cost function value of the dual problem is the same as that of the primal problem. The $\lambda$'s that correspond to the power balance equations at every bus are commonly referred to as the locational marginal price (LMP) value at that bus. These represent the price that an entity at that node would have to pay in order to consume an additional unit of power. For example, assume we will constrain the price at bus $k$ by adding the constraint $\lambda_k \leq m$ to the optimization problem. Adding a constraint into the dual problem adds an additional variable $\alpha$ into the primal problem. Assuming the power balance constraint for bus $k$ is constraint $i$, the new primal problem is thus

\vspace{-.2cm}
\begin{small}
\begin{equation} 
\begin{aligned}
\hspace{-.5cm} \underset{x}{\text{minimize}}  \hspace{.2cm}
\frac{1}{2}x^TQx + &c^Tx + m\alpha  \\
\text{subject to} \hspace{0.6cm}
(Ax)_j &\leq b_j, j \neq i \\
(Ax)_i - \alpha &\leq b_i \\
\alpha &\geq 0
\end{aligned}
\end{equation}
\end{small}

\noindent where $(Ax)_j$ denotes row $j$ of matrix $A$ multiplied by the $j$th variable in vector $x$ and $b_j$ denotes the $j$th entry in vector $b$. Of course, the inclusion of this additional constraint introduces the possibility of making the problem infeasible. Future work into market structures and bid clearing mechanisms under this framework will be performed to ensure feasibility.

\section{Linear example}
\label{sec:example}
\subsection{Optimization Formulation}
For simplicity and clarity, first we formulate an economic dispatch problem where the network constraints are not considered and the cost function is linear. Assume an industrial customer wants to determine how many widgets to produce in a given day. The cost to run their gas generator $P_G$ to produce these widgets is \$$a$/kW. The benefit they gain from producing widgets can be defined as \$$b$/kW. The minimum and maximum amount of power they can consume is $\underline{P_L}$ and $\overline{P_L}$, respectively. Thus, the overall optimization problem is:

\vspace{-.1cm}
\begin{small}
\begin{equation} \label{ex_primal}
\begin{aligned}
 \underset{P_G, P_L}{\text{maximize}} \hspace{.5cm} -aP_G + &bP_L\\
\hspace{.2cm} \text{subject to} \hspace{.5cm}
  \underline{P_L} \leq P_L &\leq  \overline{P_L} \\
 \hspace{-1cm} P_G - P_L &= 0 \\
   P_G, P_L &\geq 0
\end{aligned}
\end{equation}
\end{small}
\vspace{-.2cm}

\noindent The dual of (\ref{ex_primal}) is then written as:

\vspace{-.1cm}
\begin{small}
\begin{equation} \label{ex_dual}
\begin{aligned}
 \underset{\lambda}{\text{minimize}} \hspace{.5cm} -\underline{P_L}\lambda_1 + \overline{P_L}\lambda_2\\
\hspace{.2cm} \text{subject to} \hspace{1.6cm}
\mu \geq -a& \\
-\lambda_1 + \lambda_2 - \mu \geq b&\\
 \lambda_1, \lambda_2, \mu \geq 0&
\end{aligned}
\end{equation}
\end{small}
\vspace{-.2cm}

\noindent Where $\mu$ is the Lagrange multiplier corresponding to the power balance constraint $P_G - P_L = 0$. Assume that the maximum amount the customer wants to pay to produce widgets for that hour is $m$; i.e., $\mu \leq m$. Reformulating the new primal problem, a variable $P_m$ corresponding to the new constraint is added, and the final new primal is as follows:

\begin{small}
\begin{equation} \label{ex_newprimal}
\begin{aligned}
 \underset{P_G, P_L, P_m}{\text{maximize}} \hspace{.8cm} -aP_G + bP_L - &mP_m\\
\hspace{.4cm} \text{subject to} \hspace{1.7cm}
  \underline{P_L} \leq P_L &\leq  \overline{P_L} \\
 \hspace{-1cm} P_G - P_L - P_m &= 0 \\
   P_G, P_L, P_m &\geq 0
\end{aligned}
\end{equation}
\end{small}
\vspace{-.5cm}

\subsection{Karush-Kuhn-Tuker (KKT) Conditions}
By comparing the KKT conditions for (\ref{ex_primal}) and (\ref{ex_newprimal}), we can observe how the optimal solution of the original primal problem is changed by the inclusion of this additional dual variable constraint, and interpret the new variable, $P_m$. Define the Lagrangian functions for the original primal and new primal problem, respectively, as follows:

\begin{small}
\begin{equation*} 
\begin{aligned}
L_1(P_G, P_L, \lambda, \mu) &= -aP_G + bP_L + \lambda_1(-P_L + \underline{P_L})\hspace{2cm} \\
&+\lambda_2(P_L-\overline{P_L}) + \mu(P_G-P_L) - \lambda_3(P_G)
\end{aligned}
\end{equation*}
\end{small}
\vspace{-.4cm}

\begin{small}
\begin{equation*} 
\begin{aligned}
L_2(P_G, P_L, P_m, \lambda, \mu) = -aP_G + bP_L - mP_m + \lambda_1(-P_L + \underline{P_L})\hspace{.2cm} \\\
+ \lambda_2(P_L-\overline{P_L}) + \mu(P_G-P_L-P_m) - \lambda_3(P_G) - \lambda_4(P_m)
\end{aligned}
\end{equation*}
\end{small}

Most of the KKT conditions of both problems are identical; however, the following relevant conditions differ: \vspace{.2cm}

\noindent \textbf{Original Primal}

\vspace{-.7cm}
\begin{small}
\begin{equation*}
\begin{aligned}
\mu(P_G-P_L) = 0
\end{aligned}
\end{equation*}
\end{small}

\noindent \textbf{New Primal}

\vspace{-.6cm}
\begin{small}
\begin{equation*}
\begin{aligned}
\mu(P_G-P_L-P_m) = 0 \\
-m-\mu-\lambda_4 = 0
\end{aligned}
\end{equation*}
\end{small}

\noindent The new variable $P_m$, which must be nonnegative, can be interpreted as the amount of load reduction required to reach the specified price, $m$, if the original optimal solution of the primal problem previously resulted in an optimal price that was greater than $m$.

\section{Conclusions and Future work}
\label{sec:conclusions}
Constraining the Lagrange multipliers in the optimal power flow problem allows for consideration of problem formulations that could not be considered previously. With this framework, advanced demand side bidding techniques could be explicitly included in the optimization problem, as well as the hedging against price fluctuations that could occur due to intermittent energy sources or the rebound effect from dynamic pricing. Future work will develop more detailed market structures by utilizing this mathematical framework, including methodologies for the acceptance/rejection of bids, and provide analyses on what opportunities this framework could open for demand response programs.



\begin{thebibliography}{9}

\bibitem{rebound} W. Zhang, K. Kalsi, J. Fuller, M. Elizondo, and D. Chassin, ``Aggregate model for heterogeneous thermostatically controlled loads with demand response," \emph{IEEE PES General Meeting}, San Diego, CA, 2012.

\bibitem{duncan2011} D. S. Callaway and I. A. Hiskens, ``Achieving Controllability of Electric Loads," \emph{Proceedings of the IEEE}, vol. 99, no. 1, pp 184-199, Jan. 2011.

\bibitem{DSB}, S. Rassenti, V. Smith and B. Wilson, ``Controlling market power and price spikes in electricity networks: Demand-side bidding," \emph{Proc. Nat. Acad. Sci}, vol. 100, no. 5, pp. 2998-3003, May 2013.

\bibitem{BoVa04}, S. Boyd and L. Vandenberghe, \emph{Convex Optimization (pdf)}, Cambridge University Press, ISBN 978-0-521-83378-3. Retrieved October 3, 2011.

\end{thebibliography}

\end{document}